\documentclass[a4paper]{amsart}
\usepackage{amssymb} 
\usepackage[T1]{fontenc}
\usepackage[latin1]{inputenc}
\usepackage{amsfonts}
\usepackage{amsxtra}
\usepackage{ae}
\usepackage[all]{xy}
\usepackage{enumerate}

\include{diagram}

\newcommand*{\ket}{\rangle}
\newcommand*{\bra}{\langle}
\newcommand*{\ad}{\mathsf{ad}}

\newcommand*{\C}{\mathcal{C}}

\newcommand*{\E}{\mathcal{E}}

\renewcommand*{\H}{\mathcal{H}}

\newcommand*{\CC}{\mathcal{CC}}
\newcommand*{\CI}{\mathcal{CI}}
\renewcommand*{\P}{\mathcal{P}}

\newcommand*{\T}{\mathcal{T}}

\newcommand*{\twisted}{\boxtimes}

\newcommand*{\DD}{\mathsf{D}}

\newcommand*{\Irr}{\mathsf{Irr}}

\renewcommand*{\max}{\mathsf{f}}
\newcommand*{\red}{\mathsf{r}}

\newcommand*{\cop}{\mathsf{cop}}
\renewcommand*{\top}{\mathsf{top}}

\newcommand*{\HH}{\mathbb{H}}
\newcommand*{\KH}{\mathbb{K}}
\newcommand*{\LH}{\mathbb{L}}

\DeclareMathOperator{\coker}{coker}

\DeclareMathOperator{\Aut}{Aut}
\DeclareMathOperator{\aut}{aut}

\DeclareMathOperator{\res}{res}
\DeclareMathOperator{\ind}{ind}

\DeclareMathOperator{\id}{id}

\numberwithin{equation}{section}
\theoremstyle{change}
\newtheorem{theorem}{Theorem}[section]
\newtheorem{prop}[theorem]{Proposition}
\newtheorem{lemma}[theorem]{Lemma}

\newtheorem{definition}[theorem]{Definition}

\begin{document}

\title[Quantum $ SU(2) $]{Quantum $ SU(2) $ and the Baum-Connes conjecture}

\author{Christian Voigt}

\address{School of Mathematics and Statistics \\
         University of Glasgow \\
         University Gardens \\
         Glasgow G12 8QW \\
         United Kingdom 
}
\email{christian.voigt@glasgow.ac.uk}


\maketitle

\begin{abstract}
We review the formulation and proof of the Baum-Connes conjecture for the dual of the quantum group $ SU_q(2) $ of Woronowicz. 
As an illustration of this result we determine the $ K $-groups of quantum automorphism groups of simple matrix algebras. 
\end{abstract}

\section{Introduction}

The Baum-Connes conjecture \cite{BC}, \cite{BCH} is a far reaching conjecture about the operator $ K $-theory of locally compact groups. It has connections 
to representation theory and harmonic analysis as well as to index theory and topology. Since its original formulation by Baum and Connes about 
thirty years ago, the conjecture has been studied from various points of view and has had important impact on the development of 
noncommutative geometry \cite{Connesbook}. \\
The aim of the conjecture is to understand the relation between two $ K $-groups, one of them being of topological nature, while the 
other one involves analysis. More precisely, let $ G $ be a second countable locally compact group and let $ A $ be a separable 
$ G $-$ C^* $-algebra. The Baum-Connes conjecture with coefficients in $ A $ asserts that the assembly map
$$ 
\mu_A: K^{\top}_*(G;A) \rightarrow K_*(G \ltimes_\red A), 
$$
is an isomorphism. Here $ K^{\top}_*(G;A) $ is the topological $ K $-theory of $ G $ with coefficients in $ A $, 
and $ K_*(G \ltimes_\red A) $ denotes the $ K $-theory of the reduced crossed product $ G \ltimes_\red A $, which by definition is the 
analytical $ K $-group. The conjecture is known to hold for large classes of groups, let us mention in particular the deep work of 
Higson-Kasparov \cite{HKatmenable} on groups with the Haagerup property, and of Lafforgue \cite{Lafforguehyperbolic} on hyperbolic groups, 
respectively. Since the left hand side of the assembly map is accessible to computations this provides a conceptual approach to determine the $ K $-groups for 
a large variety of group $ C^* $-algebras and crossed products. \\
It is natural to ask what happens if the group $ G $ in the conjecture is replaced by a locally compact quantum group \cite{KVLCQG}. 
Although this question does not have direct connections to classical problems in topology or geometry, it is interesting from the point 
of view of operator $ K $-theory. Indeed, quantum groups and their crossed products give rise to a large class of $ C^* $-algebras, 
and $ K $-theory computations in this context are typically difficult. Since many considerations for groups generalize to the setting of quantum groups, 
one may hope that methods from the Baum-Connes conjecture extend to this broader context. \\
In this note we shall review some steps taken in this direction during the last years. We focus in particular on the case of 
the quantum group $ SU_q(2) $ of Woronowicz \cite{Woronowiczsuq2}, one of the most prominent examples in the theory of quantum groups. 
The Baum-Connes problem in this setting is rather a question about the dual discrete quantum group, and we shall describe 
the proof of the Baum-Connes conjecture for the dual of $ SU_q(2) $ given in \cite{Voigtbcfo}. In addition 
we shall explain basic facts concerning braided tensor products. The material covered here is mostly taken 
from \cite{NVpoincare}, \cite{Voigtbcfo}. We have added various comments, and expanded some aspects 
that have been treated only briefly in these papers. \\ 
Let us point out that already the correct definition of an assembly map for quantum groups presents a nontrivial problem.   
The definition of the left hand side of the Baum-Connes conjecture for groups given in \cite{BCH} involves the universal space for 
proper actions, a concept which does not translate to the quantum setting in an obvious way. An important insight due to Meyer and Nest \cite{MNtriangulated} 
is that a solution to this problem should be based on a categorical approach to proper actions. 
In fact, it turns out that the setup of Meyer and Nest is well-suited to study the Baum-Connes problem for the dual 
of $ SU_q(2) $. Most importantly, one obtains explicit $ K $-theory computations as a consequence. We illustrate this in the case of 
quantum automorphism groups of full matrix algebras in the sense of Wang \cite{Wangqsymmetry}. \\
In fact, what makes the strong Baum-Connes property for the dual of $ SU_q(2) $ particularly useful is that it passes to arbitrary free orthogonal 
quantum groups \cite{Wangfree} by monoidal equivalences \cite{BdRV}. The strategy of transporting structural results under monoidal equivalences 
has been successfully applied in other contexts as well, see \cite{VaesVergnioux}, \cite{dRV} and the recent 
paper \cite{Bichonyd}. \\
Let us indicate how this note is organized. In section \ref{secintegers} we give a short introduction 
to the categorical approach of Meyer and Nest by describing the Baum-Connes conjecture for the group $ \mathbb{Z} $. 
In section \ref{secdd} we discuss braided tensor products and their relation to the Drinfeld double. In particular, we shall explain 
the connection with the corresponding purely algebraic constructions for Hopf algebras. 
In section \ref{secsuq2} we review briefly the definition of $ SU_q(2) $ and the standard Podle\'s sphere $ SU_q(2)/T $. 
Moreover we discuss the crucial ingredient in the proof of the strong Baum-Connes property for the dual of $ SU_q(2) $, which 
amounts to a result on the equivariant $ KK $-theory of the Podle\'s sphere. 
Based on this we explain in section \ref{secbcsuq2} how to prove the Baum-Connes property for 
the dual of $ SU_q(2) $. Finally, as indicated above, we discuss how to compute the $ K $-groups of quantum automorphism groups 
of simple matrix algebras. \\
Throughout we shall use the notation adopted in \cite{NVpoincare} and \cite{Voigtbcfo}.

\section{The Baum-Connes conjecture for $ \mathbb{Z} $} \label{secintegers}

In this section we give a brief introduction to the Baum-Connes conjecture by discussing the case of the group $ \mathbb{Z} $ along 
the lines of the general theory developed by Meyer and Nest. This illustrates several features that show up in the Baum-Connes 
problem for the dual of $ SU_q(2) $ as well. In fact, the latter quantum group can be viewed as being freely generated by a single generator, in 
a similar way as $ \mathbb{Z} $ is the free group on one generator. Moreover the strong Baum-Connes property for $ \mathbb{Z} $ is actually used 
in the proof of the corresponding result for the dual of $ SU_q(2) $ described below. Throughout we will work in the framework of 
equivariant $ KK $-theory. For background information we refer the reader to \cite{Blackadarbook}. \\ 
Let us begin with some general facts and notation. If $ G $ is a second countable locally compact group we denote by $ KK^G $ the category 
defined as follows. The objects of $ KK^G $ are all separable $ G $-$ C^* $-algebras, that is, all separable $ C^* $-algebras equipped with 
a strongly continuous action of $ G $ by $ * $-automorphisms. The morphism set between two 
objects $ A $ and $ B $ is the equivariant Kasparov group $ KK^G(A,B) $, and composition of morphisms is given by Kasparov product. \\ 
Given a closed subgroup $ H \subset G $ there are two important functors relating $ KK^G $ and $ KK^H $. Firstly, we have the restriction 
functor $ \res^G_H: KK^G \rightarrow KK^H $ which 
is obtained by restricting the group action in the obvious way. Secondly, we have the induction functor $ \ind_H^G: KK^H \rightarrow KK^G $, which 
on the level of objects associates to an $ H $-$ C^* $-algebra $ B $ the induced $ G $-$ C^* $-algebra 
\begin{align*}
\ind_H^G(B) = \{f \in C_b(G, B)| f(xs) &= s^{-1} \cdot f(x) \; \text{for all} \; x \in G, s \in H \\
&\;\text{and}\; xH \mapsto ||f(xH)|| \in C_0(G/H) \}.
\end{align*}
The relations between these functors for various subgroups play a central role in the categorical approach to the Baum-Connes 
conjecture \cite{MNtriangulated}. More specifically, one considers the following full subcategories 
of $ KK^G $, 
\begin{align*}
\CC_G &= \{A \in KK^G|\res^G_H(A) = 0 \in KK^H \; \text{for all compact subgroups}\; H \subset G \} \\
\CI_G &= \{\ind_H^G(B)|B \in KK^H \; \text{for some compact subgroup} \; H \subset G \}, 
\end{align*}
and refers to their objects as compactly contractible and compactly induced $ G $-$ C^* $-algebras, respectively. \\
The starting point of the work of Meyer and Nest is the fact that $ KK^G $ is a triangulated category in a natural way. We shall not 
go into details here, basically, the triangulated structure consists of exact triangles which encode exact sequences, 
and a translation functor which associates to a $ G $-$ C^* $-algebra $ A $ its suspension $ \Sigma A = C_0(\mathbb{R}) \otimes A $. 
The localising subcategory $ \bra \CI_G \ket $ of $ KK^G $ generated by $ \CI_G $ plays a particularly important role. Roughly speaking, 
this is the full subcategory consisting of all objects that can be built from $ \CI_G $ by taking exact triangles, suspensions and countable direct sums. 
Meyer and Nest show that the categories $ \CC_G $ and $ \bra \CI_G \ket $ form a complementary pair \cite{MNtriangulated}. 
This is closely related to the existence of Dirac morphisms and the definition of the Baum-Connes assembly map. \\ 
Let us explain these constructions concretely in the special case of the group $ \mathbb{Z} $. Clearly, the only compact subgroup of $ \mathbb{Z} $ is 
the trivial group, and accordingly the category $ \CI_\mathbb{Z} $ consists of all $ \mathbb{Z} $-$ C^* $-algebras of the form 
$ C_0(\mathbb{Z}) \otimes B $ where $ B $ is any separable $ C^* $-algebra. Here $ B $ is viewed as a trivial $ \mathbb{Z} $-$ C^* $-algebra 
and $ \mathbb{Z} $ acts on the first tensor factor $ C_0(\mathbb{Z}) $ by translation. \\
The Dirac element for $ \mathbb{Z} $ is obtained from the Dirac operator on the real line, thus explaining the terminology used 
in the general setup. More precisely, the Dirac operator on $ \mathbb{R} $ is the self-adjoint unbounded operator acting in 
$ L^2(\mathbb{R}) $ by standard differentiation on smooth functions with compact support. This operator defines an odd equivariant 
$ K $-homology class for $ C_0(\mathbb{R}) $ equipped with the translation action of $ \mathbb{Z} $. Using suspension we may 
write this class as an element $ D \in KK^\mathbb{Z}(\P, \mathbb{C}) $ where $ \P = \Sigma C_0(\mathbb{R}) $. \\
The fact that the space $ \mathbb{R} $ shows up at this point is not a coincidence, in fact $ \mathbb{R} = \E \mathbb{Z} $ is the universal 
proper $ \mathbb{Z} $-space featuring in the usual definition of the Baum-Connes assembly map \cite{BCH}. The space $ \E \mathbb{Z} $ 
is related to $ \P $ by Poincar\'e duality \cite{Kasparov2}, and we may view the Dirac element $ D $ as a replacement of the canonical 
map $ \E \mathbb{Z} \rightarrow \star $ to the one-point space. \\
Now let $ A $ be a separable $ \mathbb{Z} $-$ C^* $-algebra and let $ D_A \in KK^{\mathbb{Z}}(\P \otimes A, A) $ be the morphism obtained by
taking the exterior product of $ D $ with the identity on $ A $. In the framework of Meyer and Nest, the assembly map for $ \mathbb{Z} $ 
with coefficients in $ A $ is the homomorphism 
$$
\mu_A: K_*(\mathbb{Z} \ltimes (\P \otimes A)) \rightarrow K_*(\mathbb{Z} \ltimes A) 
$$
induced by $ D_A $ after taking crossed products. Note that we do not have to distinguish between 
full and reduced crossed products here since the group $ \mathbb{Z} $ is amenable. \\
The Baum-Connes conjecture for $ \mathbb{Z} $ with coefficients in $ A $ asserts that $ \mu_A $ is an isomorphism. In fact, 
the following strong Baum-Connes property holds in this case. 
\begin{theorem} \label{DiracZ}
Let $ A $ be a separable $ \mathbb{Z} $-$ C^* $-algebra. Then $ D_A \in KK^\mathbb{Z}(\P \otimes A, A) $ is invertible. 
\end{theorem} 
The proof of theorem \ref{DiracZ} is a basic instance of the Dirac-dual Dirac method of Kasparov \cite{Kasparov2}. 
There exists a dual Dirac element $ \eta \in KK^\mathbb{Z}(\mathbb{C}, \P) $, and the main step of the argument consists 
in showing that $ \eta $ is the inverse of $ D $ in the category $ KK^\mathbb{Z} $. This can be viewed as an equivariant 
version of Bott periodicity, proving the claim for $ A = \mathbb{C} $. The general case follows by taking 
exterior tensor products with the algebra $ A $. \\
An equivalent, more categorical way to formulate theorem \ref{DiracZ} is to say that the localising category 
$ \bra \CI_\mathbb{Z} \ket $ is equal to $ KK^\mathbb{Z} $. Making this explicit leads to the Pimsner-Voiculescu exact sequence for 
the $ K $-theory of crossed products by $ \mathbb{Z} $. More precisely, for every $ A \in KK^{\mathbb{Z}} $ we have an extension 
$$
\xymatrix{
0 \ar@{->}[r] & \Sigma C_0(\mathbb{Z}) \otimes A \ar@{->}[r] & C_0(\mathbb{R}) \otimes A \ar@{->}[r] & C_0(\mathbb{Z}) \otimes A \ar@{->}[r] & 0 
}
$$
of $ \mathbb{Z} $-algebras induced from the inclusion $ \mathbb{Z} \subset \mathbb{R} $, here $ \Sigma $ denotes suspension as above. 
Taking crossed products with $ \mathbb{Z} $ and applying $ K $-theory yields a six-term exact sequence of $ K $-groups. 
By Takesaki-Takai duality and stability we obtain 
$$ 
K_*(\mathbb{Z} \ltimes (C_0(\mathbb{Z}) \otimes A)) \cong K_*(\KH(l^2(\mathbb{Z})) \otimes A) \cong K_*(A). 
$$
Using theorem \ref{DiracZ} we may identify $ K_*(\mathbb{Z} \ltimes (C_0(\mathbb{R}) \otimes A)) $ with 
$ K_{* + 1}(\mathbb{Z} \ltimes A) $. This yields the Pimsner-Voiculescu exact sequence 
$$
\xymatrix{
K_0(A) \ar@{->}[r]^{\!\!\! \id - \alpha_*} \ar@{<-}[d] & K_0(A) \ar@{->}[r] &
        K_0(\mathbb{Z} \ltimes A) \ar@{->}[d] \\
   K_1(\mathbb{Z} \ltimes A) \ar@{<-}[r] &
    K_1(A) \ar@{<-}[r]^{\quad \id - \alpha_*} &
   K_1(A) \\
}
$$
where $ \alpha \in \Aut(A) $ is the automorphism implementing the action of $ \mathbb{Z} $ on $ A $. \\
In categorical language, the existence of the above extension of $ \mathbb{Z} $-$ C^* $-algebras shows that 
every $ A \in KK^{\mathbb{Z}} $ has a projective resolution of length $ 1 $. Using the strong Baum-Connes property of $ \mathbb{Z} $ 
one obtains an exact triangle of the form 
$$
\xymatrix{
C_0(\mathbb{Z}) \otimes A \ar@{->}[r] & C_0(\mathbb{Z}) \otimes A \ar@{->}[r] & A \ar@{->}[r] & \Sigma C_0(\mathbb{Z}) \otimes A
}
$$
for every $ A \in KK^{\mathbb{Z}} $. \\
The basic argument leading to the Pimsner-Voiculescu sequence works in much greater generality \cite{MNhomalg1}. 
One of the main results in \cite{VVfreeu} is an analogue of the Pimnser-Voiculescu sequence for free quantum groups, 
obtained from the strong Baum-Connes property for these quantum groups in the same way as above.

\section{Braided tensor products and the Drinfeld double} \label{secdd}

Before focussing on $ SU_q(2) $ we shall discuss in this section a specific problem with coactions of quantum groups which appears naturally in 
connection with the Baum-Connes conjecture. This problem does not show up in the classical case and might at first glance be surprising. For 
the technical details we refer to \cite{NVpoincare}. \\
Let us again consider the Baum-Connes conjecture for the group 
$ \mathbb{Z} $. Recall from section \ref{secintegers} that the proof of the strong Baum-Connes property consists of two steps in this case. 
The first, and crucial part of the proof of theorem \ref{DiracZ} is to show that the Dirac element $ D \in KK^{\mathbb{Z}}(\P, \mathbb{C}) $ 
is invertible. In the second part of the proof one takes exterior products to extend this to arbitrary coefficient algebras. \\
If we want to follow a similar strategy for a quantum group $ G $, the second, rather formal step of the argument 
turns out to be problematic. In fact, there is no natural tensor product operation on the category 
of $ G $-$ C^* $-algebras in general. To circumvent this one is naturally lead to study braided tensor products and actions of the Drinfeld 
double of $ G $. In the sequel we shall explain these constructions and indicate the link with well-known results in the algebraic setting 
of Hopf algebras. \\
Indeed, the basic problem with tensor products of coactions is purely algebraic, and it can be most efficiently explained using the 
language of monoidal categories. Let $ \C $ be a monoidal category, which for simplicity we assume to be strict. By definition, an algebra 
in $ \C $ is an object $ A \in \C $ together with a morphism $ \mu_A: A \otimes A \rightarrow A $ such that the diagram 
$$
\xymatrix{
A \otimes A \otimes A \ar@{->}[r]^{\quad \id \otimes \mu_A} \ar@{->}[d]^{\mu_A \otimes \id} & A \otimes A \ar@{->}[d]^{\mu_A} \\
A \otimes A \ar@{->}[r]^{\mu_A} & A
     }
$$
is commutative. This definition amounts of course to an algebra without unit, but the existence of units does not affect our discussion. 
Assume now that $ A $ and $ B $ are algebras in $ \C $. We may form the tensor product $ A \otimes B $ as an object of $ \C $, 
but in contrast to the situation for, say the category of vector spaces over a field, this object will typically not be an algebra 
in $ \C $ in a natural way. What is needed is a prescription how to exchange the order of tensor products. \\ 
The situation changes if the monoidal category $ \C $ is braided. If we assume that $ \C $ is braided and $ \gamma_{BA}: B \otimes A \rightarrow A \otimes B $ 
denotes the braiding, a natural multiplication $ \mu_{A \otimes B} $ for $ A \otimes B $ is defined as the composition
$$
\xymatrix{
A \otimes B \otimes A \otimes B \; \; \; \ar@{->}[r]^{\id \otimes \gamma_{BA} \otimes \id} & \; \; A \otimes A \otimes B \otimes B 
\ar@{->}[r]^{\quad \quad \; \; \mu_A \otimes \mu_B} & A \otimes B.
     }
$$
For instance, this yields the usual tensor product algebra structure if $ \C $ is the category of $ G $-modules for a discrete group 
$ G $ with the braiding given by the flip map. \\
However, the monoidal categories we have to work with are usually far from being braided, 
even with the notion of braiding interpreted in a loose sense. It is therefore important that the above construction of a tensor product algebra 
still works if one of the objects is an algebra in the Drinfeld center $ Z(\C) $ of the category $ \C $. \\ 
The Drinfeld center of a monoidal category $ \C $ is a braided monoidal category whose objects are objects of $ \C $ together with a specified way 
of permuting them with arbitrary objects of $ \C $ in tensor products, see \cite{JStortileYB}, \cite{Kasselbook}. In the case 
that $ \C $ is the category of modules over a Hopf algebra $ H $, the Drinfeld center of $ \C $ is the category of $ H $-Yetter-Drinfeld modules. 
If in addition $ H $ is finite dimensional, 
the latter is equivalent to the category of modules over the Drinfeld double of $ H $. \\
We are interested in a situation where, loosely speaking, the Hopf algebra $ H $ is replaced by a locally compact quantum group. In this 
generality the above picture has to be adapted appropriately, but it should be no surprise that this leads to Yetter-Drinfeld structures 
and the Drinfeld double in the operator algebraic framework. \\ 
To explain the analogy with the algebraic theory let us recall some definitions. If $ H $ is a Hopf algebra we use the Sweedler 
notation $ \Delta(x) = x_{(1)} \otimes x_{(2)} $ for the  comultiplication. Moreover we write $ S $ and $ \epsilon $ for the antipode and 
counit of $ H $, respectively. For the sake of definiteness we shall work over the complex numbers. 
\begin{definition} \label{defydalg}
Let $ H $ be a Hopf algebra. An $ H $-Yetter-Drinfeld module is a vector space $ M $ 
which is both a left $ H $-module via $ H \otimes M \rightarrow M, f \otimes m \mapsto f \cdot m $ and 
a left $ H $-comodule via $ M \rightarrow H \otimes M, m \mapsto m_{(-1)} \otimes m_{(0)} $ such that 
$$
(f \cdot m)_{(-1)} \otimes (f \cdot m)_{(0)} = f_{(1)} m_{(-1)} S(f_{(3)}) \otimes f_{(2)} \cdot m_{(0)} 
$$
for all $ f \in H $ and $ m \in M $. 
\end{definition}
If $ H $ is a finite dimensional Hopf algebra we write $ H^* $ for the dual Hopf algebra. Moreover 
we let $ (H^*)^\cop $ be the Hopf algebra obtained by considering $ H^* $ with the opposite coproduct. We shall write 
$ \hat{\Delta} $ for the coproduct of $ (H^*)^\cop $. Let $ e_1, \dots, e_n $ be a basis of $ H $ with dual 
basis $ e^1, \dots, e^n $ of $ (H^*)^\cop = H^* $ and consider the element 
$$
w = \sum_{j = 1}^n e_j \otimes e^j \in H \otimes (H^*)^\cop.  
$$
The following properties of $ w $ can be verified by direct calculation. 
\begin{lemma} \label{bichar}
The element $ w $ is a bicharacter of $ H \otimes (H^*)^\cop $, that is, $ w $ is invertible and the formulas 
$$
(\epsilon \otimes \id)(w) = 1, \qquad (\id \otimes \epsilon)(w) = 1 
$$
as well as 
$$
(\Delta \otimes \id)(w) = w_{13} w_{23}, \qquad (\id \otimes \hat{\Delta})(w) = w_{13} w_{12} 
$$
hold. 
\end{lemma}
Here we have used the leg numbering notation. The definition of a Yetter-Drinfeld module can now be rephrased as follows. 
\begin{lemma} \label{ydcomod}
Let $ H $ be a finite dimensional Hopf algebra. Then an $ H $-Yetter-Drinfeld module is the same thing as 
a vector space $ M $ which is both a left $ H $-comodule via $ \alpha: M \rightarrow H \otimes M $ 
and a left $ (H^*)^\cop $-comodule via $ \lambda: M \rightarrow (H^*)^\cop \otimes M $ such that the diagram
$$
\xymatrix{
M \ar@{->}[r]^{\!\!\!\!\!\!\!\!\! \lambda} \ar@{->}[d]^\alpha & (H^*)^\cop \otimes M \ar@{->}[r]^{\!\!\!\!\!\!\!\!\! \id \otimes \alpha} &
(H^*)^\cop \otimes H \otimes M \ar@{->}[d]^{\sigma \otimes \id} \\
H \otimes M \ar@{->}[r]^{\!\!\!\!\!\!\!\!\!\!\!\!\! \id \otimes \lambda} & H \otimes (H^*)^\cop \otimes M\; \ar@{->}[r]^{\ad(w) \otimes \id}
& \; H \otimes (H^*)^\cop \otimes M
     }
$$
is commutative.
\end{lemma} 
The correspondence is given by identifying the coaction $ \lambda: M \rightarrow (H^*)^\cop \otimes M  $ with a left $ H $-module structure 
on $ M $ by duality. 
\begin{definition} 
Let $ H $ be a finite dimensional Hopf algebra. The Drinfeld codouble of $ H $ is 
$$
\DD_H = H \otimes (H^*)^\cop
$$
with the tensor product algebra structure, the comultiplication 
\begin{align*}
\Delta_\DD(f \otimes x) &= (\id \otimes \sigma \otimes \id)(\id \otimes \ad(w) \otimes \id)(\Delta \otimes \hat{\Delta})(f \otimes x),  
\end{align*} 
the counit $ \epsilon_\DD(f \otimes x) = \epsilon(f) \epsilon(x) $ and the antipode $ S_\DD = (S \otimes S) \ad(w) $, 
where $ \ad(w) $ denotes conjugation by $ w $ and $ \sigma $ is the flip map. 
\end{definition} 
Using lemma \ref{bichar} it is straightforward to check that $ \DD_H $ is a Hopf algebra such that the canonical 
projection maps $ \pi: \DD_H \rightarrow H, \pi(f \otimes x) = f \epsilon(x) $ and 
$ \hat{\pi}: \DD_H \rightarrow (H^*)^\cop, \hat{\pi}(f \otimes x) = \epsilon(f) x $ are Hopf algebra homomorphisms. 
The dual Hopf algebra of the Drinfeld codouble is the Drinfeld double of $ H $, in the algebraic context 
the double is usually studied from this dual point of view. \\
The following basic result explains the connection between Yetter-Drinfeld modules and comodules over the Drinfeld double, see
for instance \cite{KS}. 
\begin{prop} \label{yddouble}
Let $ H $ be a finite dimensional Hopf algebra. Then there is a bijective correspondence between $ H $-Yetter-Drinfeld modules and 
left $ \DD_H $-comodules. 
\end{prop} 
The correspondence is given as follows. If $ \lambda: M \rightarrow \DD_H \otimes M $ is a $ \DD_H $-comodule structure, then the 
Hopf algebra homomorphisms $ \pi $ and $ \hat{\pi} $ defined above induce coactions $ M \rightarrow H \otimes M $ and 
$ M \rightarrow (H^*)^\cop \otimes M $. These coactions satisfy the compatibility relation in lemma \ref{ydcomod}. \\
From proposition \ref{yddouble} it follows in particular that the category of $ H $-Yetter-Drinfeld modules over a 
finite dimensional Hopf algebra $ H $ is a monoidal category in a 
natural way. Note that an algebra in this category can be defined as an algebra $ A $ which is both an $ H $-comodule 
algebra and a $ (H^*)^\cop $-comodule algebra such that the compatibility condition in lemma \ref{ydcomod} holds for $ M = A $. \\
Let us now go back to $ C^* $-algebras. If $ G $ is a locally compact quantum group we write 
$ W \in M(C^\red_0(G) \otimes C^*_\red(G)) $ for the fundamental multiplicative unitary \cite{KVLCQG}. This unitary replaces the 
element $ w $ considered in the algebraic setting above. \\
The analogue of definition \ref{defydalg} for actions on $ C^* $-algebras reads as follows. 
\begin{definition} \label{defydcstar}
Let $ G $ be a locally compact quantum group and let $ S = C^{\red}_0(G) $ and $ \hat{S} = C^*_\red(G) $ be the
associated reduced Hopf-$ C^* $-algebras. A $ G $-Yetter-Drinfeld $ C^* $-algebra is a
$ C^* $-algebra $ A $ equipped with continuous coactions $ \alpha $ of $ S $ and $ \lambda $ of $ \hat{S} $ such that the diagram
$$
\xymatrix{
A \ar@{->}[r]^\lambda \ar@{->}[d]^\alpha & M(\hat{S} \otimes A) \ar@{->}[r]^{\!\!\!\!\!\!\!\!\! \id \otimes \alpha} &
M(\hat{S} \otimes S \otimes A) \ar@{->}[d]^{\sigma \otimes \id} \\
M(S \otimes A) \ar@{->}[r]^{\!\!\!\!\!\!\!\!\!\! \id \otimes \lambda} & M(S \otimes \hat{S} \otimes A)\; \ar@{->}[r]^{\ad(W) \otimes \id}
& \;M(S \otimes \hat{S} \otimes A)
     }
$$
is commutative. 
\end{definition}
We may define Yetter-Drinfeld actions on Hilbert spaces or Hilbert modules in a similar way, thus obtaining an even closer analogy 
to the algebraic constructions above. \\
If $ G $ is a locally compact quantum group, then the Drinfeld double $ \DD(G) $ of $ G $ is given by 
$ C_0^\red(\DD(G)) = C_0^\red(G) \otimes C^*_\red(G) $ with the comultiplication
$$
\Delta_{\DD(G)} = (\id \otimes \sigma \otimes \id)(\id \otimes \ad(W) \otimes \id)(\Delta \otimes \hat{\Delta})
$$
where $ \ad(W) $ denotes the adjoint action of $ W $ and $ \sigma $ is the flip map. Comparing this with the 
algebraic setting one should keep in mind that in the conventions of Kustermans and Vaes the comultiplication of 
$ C^*_\red(G) $ is already flipped by default. \\
We have the following analogue of proposition \ref{yddouble}, see \cite{NVpoincare}. 
\begin{prop} \label{yddoublecstar}
Let $ G $ be a locally compact quantum group and let $ \DD(G) $ be its Drinfeld double. Then a $ G $-Yetter-Drinfeld $ C^* $-algebra is
the same thing as a $ \DD(G) $-$ C^* $-algebra.
\end{prop}
We shall now define the braided tensor product of a $ G $-Yetter-Drinfeld $ C^* $-algebra $ A $ with a $ G $-$ C^* $-algebra $ B $. 
For this construction it is in fact not necessary to write down the braiding of the Drinfeld double.
Let $ \HH = L^2(G) $ be the GNS-space of the left Haar weight of $ G $, so that 
$ C^\red_0(G) $ and $ C^*_\red(G) $ are naturally $ C^* $-subalgebras of $ \LH(\HH) $. If $ \beta: B \rightarrow M(C^\red_0(G) \otimes B) $ 
implements the action of $ G $ then $ B $ acts on the Hilbert module $ \HH \otimes B $ by $ \beta $. Similarly, if 
$ \lambda: A \rightarrow M(C^*_\red(G) \otimes A) $ is the coaction of $ C^*_\red(G) $ on $ A $ then $ A $ acts on $ \HH \otimes A $ by
$ \lambda $. 
From this we obtain two $ * $-homomorphisms $ \iota_A = \lambda_{12}: A \rightarrow \LH(\HH \otimes A \otimes B) $ and
$ \iota_B = \beta_{13}: B \rightarrow \LH(\HH \otimes A \otimes B) $.
\begin{definition} \label{defbraidedtensor}
Let $ G $ be a locally compact quantum group, let $ A $ be a $ G $-Yetter-Drinfeld-$ C^* $-algebra and $ B $ a $ G $-$ C^* $-algebra.
With the notation as above, the braided tensor product $ A \twisted_G B $ is the
$ C^* $-subalgebra of $ \LH(\HH \otimes A \otimes B) $ generated by all elements
$ \iota_A(a) \iota_B(b) $ for $ a \in A $ and $ b \in B $.
\end{definition}
It turns out that the braided tensor product $ A \twisted_G B $ is in fact equal to the closed linear span $ [\iota_A(A) \iota_B(B)] $. 
In particular, we have natural nondegenerate $ * $-homomorphisms
$ \iota_A: A \rightarrow M(A \twisted B) $ and $ \iota_B: B \rightarrow M(A \twisted B) $. \\
The braided tensor product shares the basic properties that hold in the algebraic setting. 
For instance, $ A \twisted B $ is a $ G $-$ C^* $-algebra in a canonical way such that the $ * $-homomorphisms
$ \iota_A $ and $ \iota_B $ are $ G $-equivariant. If $ B $ is a $ \DD(G) $-$ C^* $-algebra then $ A \twisted B $ is 
a $ \DD(G) $-$ C^* $-algebra such that $ \iota_A $ and $ \iota_B $ are $ \DD(G) $-equivariant. \\
Observe that the braided tensor product defined above generalizes the minimal tensor product of $ C^* $-algebras. 
We may refer to it as the minimal braided tensor product.

\section{The quantum group $ SU_q(2) $ and the Podle\'s sphere} \label{secsuq2}

In this section we discuss some results related to $ SU_q(2) $ and the standard Podle\'s sphere that constitute 
the core of the proof of the Baum-Connes conjecture. For simplicity we shall restrict to the case $ q \in (0,1] $ in the sequel, 
although the main arguments work with minor modifications for negative deformation parameters as well. For background material on 
quantum groups we refer to \cite{Kasselbook}, \cite{KS}. \\
Let us first recall the definition of $ SU_q(2) $, see \cite{Woronowiczsuq2}. 
\begin{definition} \label{defsuq2}
The $ C^* $-algebra $ C(SU_q(2)) $ is the universal $ C^* $-algebra generated by elements 
$ \alpha $ and $ \gamma $ satisfying the relations
$$
\alpha \gamma = q \gamma \alpha, \quad \alpha \gamma^* = q \gamma^* \alpha, \quad \gamma \gamma^* = \gamma^* \gamma, \quad 
\alpha^* \alpha + \gamma^* \gamma = 1, \quad \alpha \alpha^* + q^2 \gamma \gamma^* = 1. 
$$
The comultiplication $ \Delta: C(SU_q(2)) \rightarrow C(SU_q(2)) \otimes C(SU_q(2)) $ is defined by 
$$
\Delta(\alpha) = \alpha \otimes \alpha - q \gamma^* \otimes \gamma, \qquad \Delta(\gamma) = \gamma \otimes \alpha + \alpha^* \otimes \gamma. 
$$
\end{definition} 
The relations in definition \ref{defsuq2} are equivalent to saying that the fundamental matrix 
$$
u = 
\begin{pmatrix}
\alpha & -q \gamma^* \\
\gamma & \alpha^*
\end{pmatrix}
$$
is unitary. \\
We write $ \mathbb{C}[SU_q(2)] $ for the Hopf-$ * $-algebra of polynomial functions on $ SU_q(2) $. By definition, this is the 
dense $ * $-subalgebra of $ C(SU_q(2)) $ generated by $ \alpha $ and $ \gamma $. 
We use Sweedler notation $ \Delta(x) = x_{(1)} \otimes x_{(2)} $ for the comultiplication, and write $ \epsilon $ and $ S $ for the 
counit and the antipode of $ \mathbb{C}[SU_q(2)] $, respectively. \\
The Hilbert space $ L^2(SU_q(2)) $ is the completion of $ C(SU_q(2)) $ with respect to the inner product 
$$ 
\bra x, y \ket = \phi(x^* y) 
$$ 
induced by the Haar state $ \phi $. It is a $ SU_q(2) $-Hilbert space with the left regular representation. 
We may choose an orthonormal basis of $ L^2(SU_q(2)) $ according to the decomposition into isotypical components. Explicitly, we 
have basis vectors $ e^{(l)}_{i,j} $ where $ l \in \tfrac{1}{2} \mathbb{N} $ and $ -l \leq i,j \leq l $ run over 
integral or half-integral values, respectively. \\
The classical torus $ T = S^1 $ is a closed quantum subgroup of $ SU_q(2) $ determined by the $ * $-homomorphism
$ \pi: \mathbb{C}[SU_q(2)] \rightarrow \mathbb{C}[T] = \mathbb{C}[z,z^{-1}] $ given in matrix notation by
$$
\pi\begin{pmatrix}
\alpha & -q \gamma^* \\
\gamma & \alpha^*
\end{pmatrix}
= 
\begin{pmatrix}
z & 0 \\
0 & z^{-1}
\end{pmatrix}. 
$$
By definition, the standard Podle\'s sphere $ SU_q(2)/T  $ is the corresponding homogeneous space \cite{Podlesspheres}. The algebra of polynomial 
functions on $ SU_q(2)/T $ is given by 
$$ 
\mathbb{C}[SU_q(2)/T] = \{x \in \mathbb{C}[SU_q(2)] | (\id \otimes \pi)\Delta(x) = x \otimes 1 \}, 
$$ 
and the $ C^* $-algebra $ C(SU_q(2)/T) $ is the closure of $ \mathbb{C}[SU_q(2)/T] $ inside $ C(SU_q(2)) $. \\
More generally, for $ k \in \mathbb{Z} $ we define 
\begin{equation*}
\Gamma(\E_k) = \{x \in \mathbb{C}[SU_q(2)] | (\id \otimes \pi)\Delta(x) = x \otimes z^k \} \subset \mathbb{C}[SU_q(2)] 
\end{equation*}
and let $ C(\E_k) $ and $ L^2(\E_k) $ be the closures of $ \Gamma(\E_k) $ in $ C(SU_q(2)) $ and $ L^2(SU_q(2)) $, respectively. 
Note that we have $ \Gamma(\E_0) = \mathbb{C}[SU_q(2)/T] $ by construction. The space $ \Gamma(\E_k) $ is a $ \mathbb{C}[SU_q(2)/T] $-bimodule in a natural way for all $ k \in \mathbb{Z} $. Using Hopf-Galois theory it can be shown that $ \Gamma(\E_k) $ is finitely generated and projective 
both as a left and right $ \mathbb{C}[SU_q(2)/T] $-module, compare \cite{Shopfgalois}. 
The space $ C(\E_k) $ is naturally a $ SU_q(2) $-equivariant Hilbert $ C(SU_q(2)/T) $-module, and $ L^2(\E_k) $ is naturally 
a $ SU_q(2) $-Hilbert space. These structures are induced from $ C(SU_q(2)) $ and $ L^2(SU_q(2)) $, respectively. \\
The above spaces admit canonical actions of the Drinfeld double $ \DD(SU_q(2)) $ of $ SU_q(2) $. 
We refer to section \ref{secdd} for the definition of the Drinfeld double and the description of its actions. 
The $ C^* $-algebra  $ C(SU_q(2)/T) $ is a $ \DD(SU_q(2)) $-$ C^* $-algebra with the action of $ SU_q(2) $ 
by translations and the coaction $ \lambda: C(SU_q(2)/T) \rightarrow M(C^*(SU_q(2)) \otimes C(SU_q(2)/T)) $ given by 
$$
\lambda(g) = \hat{W}^*(1 \otimes g) \hat{W}.
$$
Here $ \hat{W} = \Sigma W^* \Sigma $ where $ W \in M(C(SU_q(2)) \otimes C^*(SU_q(2))) $ is the fundamental multiplicative unitary and 
$ \Sigma $ is the flip map. 
The coaction $ \lambda $ is determined on the algebraic level by the adjoint action 
$$
h \cdot g = h_{(1)} g S(h_{(2)}) 
$$
of $ \mathbb{C}[SU_q(2)] $ on $ \mathbb{C}[SU_q(2)/T] $. 
The same construction turns the spaces $ C(\E_k) $ for $ k \in \mathbb{Z} $ into $ \DD(SU_q(2)) $-equivariant Hilbert $ C(SU_q(2)/T) $-modules 
for every $ k \in \mathbb{Z} $. In the case of the Hilbert spaces $ L^2(\E_k) $ we have to twist the above formula to take into account the 
fact that the Haar state $ \phi $ is not a trace in general, see \cite{Voigtbcfo}. \\
Our aim is to describe the Podle\'s sphere as an element in the equivariant $ KK $-category $ KK^{\DD(SU_q(2))} $. 
It is well-known that the $ C^* $-algebra $ C(SU_q(2)/T) $ of the Podle\'s sphere is isomorphic to $ \KH^+ $ for $ q \neq 1 $,  
the algebra $ \KH $ of compact operators on a separable Hilbert space with a unit adjoined. Using this fact it is easy to show 
that $ C(SU_q(2)/T) $ is isomorphic to $ \mathbb{C} \oplus \mathbb{C} $ in the category $ KK $. 
However, the most obvious such isomorphism does not respect the $ \DD(SU_q(2)) $-actions, in fact not even the canonical $ SU_q(2) $-actions 
on both sides. In order to obtain the desired statement on the level of $ KK^{\DD(SU_q(2))} $ we need more refined arguments. \\
More precisely, we have to work with the equivariant Fredholm module corresponding to the Dirac operator on the standard 
Podle\'s sphere, compare \cite{DSpodles}, \cite{NVpoincare}. The underlying graded $ SU_q(2) $-Hilbert space is 
$$ 
\H = L^2(\E_1) \oplus L^2(\E_{-1})
$$ 
as defined above. The representation $ \mu $ of $ C(SU_q(2)/T) $ is given by left multiplication. 
We obtain a $ G $-equivariant self-adjoint unitary operator $ F $ on $ \H $ by 
$$
F =
\begin{pmatrix}
0 & 1 \\
1 & 0
\end{pmatrix}
$$
by identifying the basis vectors $ e^{(l)}_{i,1/2} $ and $ e^{(l)}_{i, -1/2} $ in even and odd degrees. Note moreover 
that the Drinfeld double $ \DD(SU_q(2)) $ acts on $ C(SU_q(2)/T) $ and $ \H $ in the way explained above.  
\begin{prop} \label{Diracconst}
The triple $ D = (\H, \mu, F) $ is a $ \DD(SU_q(2)) $-equivariant Fredholm module defining an element $ [D] $ in 
$ KK^{\DD(SU_q(2))}(C(SU_q(2)/T), \mathbb{C}) $. 
\end{prop} 
We have already mentioned above that $ C(\E_k) $ is a $ \DD(SU_q(2)) $-equivariant Hilbert $ C(SU_q(2)/T) $-module
in a natural way. Left multiplication yields a $ \DD(SU_q(2)) $-equivariant $ * $-homomorphism $ \psi: C(SU_q(2)/T) \rightarrow \KH(C(\E_k)) $, 
and it is easy to check that $ (C(\E_k), \psi, 0) $ defines a class $ [[\E_k]] $ in
$ KK^{\DD(SU_q(2))}(C(SU_q(2)/T), C(SU_q(2)/T)) $. Moreover, for the Kasparov product of these elements the relation 
$$ 
[[\E_m]] \otimes_{C(SU_q(2)/T)} [[\E_n]] = [[\E_{m + n}]] 
$$ 
holds for all $ m,n \in \mathbb{Z} $. \\
Let us now define classes $ [D_k] \in KK^{\DD(SU_q(2))}(C(SU_q(2)/T), \mathbb{C}) $ corresponding to twisted Dirac operators 
on $ SU_q(2)/T $. More precisely, we consider the Kasparov product 
$$
[D_k] = [[\E_k]] \otimes_{C(SU_q(2)/T)} [D]
$$
where $ [D] \in KK^{\DD(SU_q(2))}(C(SU_q(2)/T), \mathbb{C}) $ is the element from proposition \ref{Diracconst}. Remark that $ [D_0] = [D] $ 
since $ [[\E_0]] = 1 $. \\
The unit homomorphism $ u: \mathbb{C} \rightarrow C(SU_q(2)/T) $ given by $ u(1) = 1 $ 
induces a class $ [u] $ in $ KK^{\DD(SU_q(2))}(\mathbb{C}, C(SU_q(2)/T)) $.
We define $ [\E_k] $ in $ KK^{\DD(SU_q(2))}(\mathbb{C}, C(SU_q(2)/T)) $ to be the Kasparov product 
$$
[\E_k] = [u] \otimes_{C(SU_q(2)/T)} [[\E_k]].
$$
Moreover, we let $ \alpha_q \in KK^{\DD(SU_q(2))}(C(SU_q(2)/T), \mathbb{C} \oplus \mathbb{C}) $ and $ \beta_q 
\in KK^{\DD(SU_q(2))}(\mathbb{C} \oplus \mathbb{C}, C(SU_q(2)/T)) $ be given by 
$$
\alpha_q = [D_0] \oplus [D_{-1}], \qquad \beta_q = (-[\E_1]) \oplus [\E_0],
$$
respectively. 
\begin{theorem} \label{PDPodles}
Let $ q \in (0,1] $. The standard Podle\'s sphere $ C(SU_q(2)/T) $ is isomorphic to 
$ \mathbb{C} \oplus \mathbb{C} $ in $ KK^{\DD(SU_q(2))} $. 
\end{theorem} 
\proof We claim that $ \beta_q $ and $ \alpha_q $ define inverse isomorphisms. 
The crucial part of the argument is the relation 
$$ 
\beta_q \circ \alpha_q = \id 
$$ 
in $ KK^{\DD(SU_q(2))}(\mathbb{C} \oplus \mathbb{C}, \mathbb{C} \oplus \mathbb{C}) $. 
In order to prove it we have to compute the Kasparov products $ [\E_0] \circ [D] $ and $ [\E_{\pm 1}] \circ [D] $. \\
The class $ [\E_0] \circ [D] $ is obtained from the $ \DD(SU_q(2)) $-equivariant Fredholm module 
$ D $ by forgetting the left action of $ C(SU_q(2)/T) $. The operator $ F $ 
intertwines the representations of $ C(SU_q(2)) $ on $ L^2(E_1) $ and $ L^2(E_{-1}) $ induced from the $ \DD(SU_q(2)) $-Hilbert space structure.
It follows that the resulting $ \DD(SU_q(2)) $-equivariant Kasparov $ \mathbb{C} $-$ \mathbb{C} $-module is degenerate, 
and hence $ [\E_0] \circ [D] = 0 $ in $ KK^{\DD(SU_q(2))}(\mathbb{C}, \mathbb{C}) $. \\
It remains to calculate $ [\E_{\pm 1}] \circ [D] $. Using $ SU_q(2) $-equivariance it is straightforward to 
show that $ [\E_{-1}] \circ [D] = 1 $ in $ KK^{SU_q(2)}(\mathbb{C} \oplus \mathbb{C}, \mathbb{C} \oplus \mathbb{C}) $. 
The entire difficulty lies in constructing a $ \DD(SU_q(2)) $-equivariant homotopy to obtain the 
same relation on the level of $ KK^{\DD(SU_q(2))} $. This can be done using explicit estimates involving Clebsch-Gordan
coefficients. For the details we refer to \cite{Voigtbcfo}. \qed \\ 
It would be nice to find a proof of theorem \ref{PDPodles} taking care of the action of the Drinfeld double 
in a more conceptual way, perhaps from a categorical point of view. Such an alternative proof might shed some light on the analogous problem in higher rank. \\
Note that for $ q = 1 $ the main difficulties in the proof of theorem \ref{PDPodles} disappear since the discrete part of the 
Drinfeld double acts trivially in this case. This is the reason why the Baum-Connes property for the dual of the classical 
group $ SU(2) $ is easier to establish than for its $ q $-deformations.

\section{The Baum-Connes conjecture for $ SU_q(2) $} \label{secbcsuq2}

In this section we discuss the proof of the Baum-Connes conjecture for the dual of the quantum 
group $ SU_q(2) $. The details of the argument can be found in \cite{Voigtbcfo}, and as in the previous section we 
shall restrict ourselves to the case $ q \in (0,1] $ for this. In the last part we explain how to compute the $ K $-groups of 
quantum automorphism groups of simple matrix algebras. \\
As discussed in section \ref{secintegers}, the formulation of the Baum-Connes conjecture in the approach 
of Meyer and Nest is based on the study of the categories of compactly contractible
and compactly induced algebras, respectively. This becomes particularly simple when there are no nontrivial compact subgroups. 
For a discrete group $ G $ this means of course that $ G $ is torsion-free. \\
It turns out that the dual of $ SU_q(2) $ is torsion-free in a suitable sense \cite{Meyerhomalg2}, \cite{Voigtbcfo}, so that we are 
in a situation which is analogous to the case of the group $ \mathbb{Z} $ explained in section \ref{secintegers}. Instead of working with the dual of 
$ SU_q(2) $ it is most convenient to use Baaj-Skandalis duality to transport the Baum-Connes problem to the 
compact side. More precisely, let us write $ G = SU_q(2) $ and, by slight abuse of notation, let us denote by $ \hat{G} $ the 
discrete quantum group determined by $ C^*(SU_q(2))^\cop = C_0(\hat{G}) $. Note that this amounts to switching the comultiplication in the 
conventions of Kustermans and Vaes. This modification is convenient for Baaj-Skandalis duality, and it 
should not lead to confusion. \\
The restriction functor from $ \hat{G} $ to the trivial quantum subgroup corresponds to 
the crossed product functor $ KK^G \rightarrow KK $ which maps $ A $ to $ G \ltimes A $ on the level of objects. Similarly, the 
induction functor from the trivial group to $ \hat{G} $ identifies with the functor 
$ \tau: KK \rightarrow KK^G $ which maps a $ C^* $-algebra $ A $ to $ \tau(A) $, the $ G $-$ C^* $-algebra obtained by 
considering the trivial action of $ G $ on $ A $. We have the following full subcategories of $ KK^G $, 
\begin{align*}
\C_G &= \{A \in KK^G| G \ltimes A = 0 \in KK \} \\
\T_G &= \{\tau(A)| A \in KK \}, 
\end{align*}
these categories correspond precisely to the compactly contractible and the compactly induced 
$ \hat{G} $-$ C^* $-algebras, respectively. These categories form a complementary pair of localising 
subcategories \cite{Meyerhomalg2}, and one can study the assembly map and the Baum-Connes problem 
for $ \hat{G} $ as for the group $ \mathbb{Z} $ in section \ref{secintegers}. 
\begin{theorem} \label{BCsuq2}
The discrete quantum group dual to $ G = SU_q(2) $ satisfies the strong Baum-Connes conjecture, that is, we have $ \bra \T_G \ket = KK^G $. 
\end{theorem}
\proof Let $ A $ be a $ G $-$ C^* $-algebra. Theorem \ref{PDPodles} implies that $ A $ is a retract of $ C(G/T) \twisted_G A $ in $ KK^G $, 
and according to the compatibility of induction with braided tensor products \cite{NVpoincare} we have a $ G $-equivariant isomorphism 
$$ 
C(G/T) \twisted_G A = \ind_T^G(\mathbb{C}) \twisted_G A \cong \ind_T^G \res_T^G(A). 
$$
As discussed in section \ref{secintegers}, the group $ \hat{T} = \mathbb{Z} $ satisfies the strong Baum-Connes conjecture. That is, we have 
$$
KK^\mathbb{Z} = \bra \CI_\mathbb{Z} \ket 
$$ 
where $ \CI_{\mathbb{Z}} $ is the full subcategory in $ KK^\mathbb{Z} $ of compactly induced $ \mathbb{Z} $-$ C^* $-algebras. 
Equivalently, we have 
$$
KK^T = \bra \T_T \ket 
$$
where $ \T_T \subset KK^T $ is the full subcategory of trivial $ T $-$ C^* $-algebras. In particular we obtain
\begin{align*} 
\res^G_T(A) \in \bra \T_T \ket \subset KK^T.  
\end{align*}
Due to theorem \ref{PDPodles} we know that 
$$ 
\ind_T^G(B) = C(G/T) \otimes B \cong (\mathbb{C} \oplus \mathbb{C}) \otimes B 
$$ 
is contained in $ \bra \T_G \ket $ inside $ KK^G $ for any trivial $ T $-$ C^* $-algebra $ B $. Since 
the induction functor $ \ind_T^G: KK^T \rightarrow KK^G $ is triangulated it therefore maps $ \bra \T_T \ket $ to $ \bra \T_G \ket $. This yields
$$
\ind_T^G \res^G_T(A) \in \bra \T_G \ket
$$
in $ KK^G $. Combining the above considerations shows $ A \in \bra \T_G \ket $, and we conclude $ KK^G = \bra \T_G \ket $
as desired. \qed \\
We remark that the case $ q = 1 $ of the previous theorem is a special case of the results in \cite{MNcompact}. \\
As already mentioned above, theorem \ref{BCsuq2} can be applied to compute the $ K $-theory of free orthogonal quantum groups \cite{Voigtbcfo}. 
If $ G $ is a free orthogonal quantum group, then the main tool for this computation is an exact triangle in $ KK^G $ of the form 
$$
\xymatrix{
C_0(G) \ar@{->}[r] & C_0(G) \ar@{->}[r] & \mathbb{C} \ar@{->}[r] & \Sigma C_0(G)
}
$$
which is analogous to the extension for the source of the Dirac morphism for $ \mathbb{Z} $ discussed in 
section \ref{secintegers}. \\
In the remaining part of this section we shall briefly discuss a further consequence of theorem \ref{BCsuq2} which is not stated in \cite{Voigtbcfo}. Let us 
consider the quantum automorphism group of $ M_n(\mathbb{C}) $ defined by Wang \cite{Wangqsymmetry}. 
By definition, this quantum group is given by the universal $ C^* $-algebra $ A_{\aut}(M_n(\mathbb{C})) $ generated by elements $ u_{ij}^{kl} $ 
for $ 1 \leq i,j,k,l \leq n $ such that 
$$
\sum_{p = 1}^n u_{ij}^{kp} u_{rs}^{pl} = \delta_{jr} u^{kl}_{is}, \qquad \sum_{p = 1}^n u_{lp}^{sr} u_{pk}^{ji} = \delta_{jr} u^{si}_{lk}
$$
and 
$$
(u_{ij}^{kl})^* = u_{ji}^{lk}, \qquad \sum_{p = 1}^n u_{kl}^{pp} = \delta_{kl}, \qquad \sum_{p = 1}^n u^{kl}_{pp} = \delta_{kl}
$$
for all $ 1 \leq i,j,k,l,r,s \leq n $. These relations are equivalent to saying that $ A_{\aut}(M_n(\mathbb{C})) $ defines a quantum group such 
that the formula
$$
\lambda(e_{ij}) = \sum_{k,l = 1}^n u^{kl}_{ij} \otimes e_{kl}
$$
determines a coaction $ \lambda: M_n(\mathbb{C}) \rightarrow A_{\aut}(M_n(\mathbb{C})) \otimes M_n(\mathbb{C}) $ 
which is trace preserving in the sense that $ (\id \otimes \tau)\lambda(x) = \tau(x) 1 $ for all $ x \in M_n(\mathbb{C}) $. 
Here $ e_{ij} $ for $ 1 \leq i,j \leq n $ are the matrix units in $ M_n(\mathbb{C}) $ and $ \tau: M_n(\mathbb{C}) \rightarrow \mathbb{C} $ 
is the standard trace. \\
Following the conventions in \cite{Voigtbcfo} we write $ A_{\aut}(M_n(\mathbb{C})) = C^*_\max(\mathbb{F} \Aut(M_n(\mathbb{C}))) $ and view 
this $ C^* $-algebra as the full group $ C^* $-algebra of a discrete quantum group $ \mathbb{F} \Aut(M_n(\mathbb{C})) $ in the sequel. 
\begin{theorem} \label{qautktheory}
Let $ n > 2 $. The discrete quantum group $ H = \mathbb{F} \Aut(M_n(\mathbb{C})) $ is $ K $-amenable and its $ K $-theory 
is given by 
$$
K_0(C^*_\max(H)) = \mathbb{Z} \oplus \mathbb{Z}_n, \qquad 
K_1(C^*_\max(H)) = \mathbb{Z},
$$
where $ \mathbb{Z}_n $ denotes the cyclic group of order $ n $. 
\end{theorem} 
\proof Let us abbreviate $ H = \mathbb{F} \Aut(M_n(\mathbb{C})) $ and write $ G = \mathbb{F}O(n) $ for the free orthogonal 
quantum group of Wang \cite{Wangfree}, see \cite{Voigtbcfo}. As mentioned in the introduction, the strong Baum-Connes property for the dual of $ SU_q(2) $ 
implies that $ G $ satisfies the strong Baum-Connes property as well. Moreover, a result of Banica \cite{Banicageneric} shows that $ H $ can be 
identified with the quantum subgroup of $ G $ generated by the coefficients of the tensor square of the fundamental corepresentation of $ G $ in the same way 
as $ C(SO(3)) $ is obtained from $ C(SU(2)) $. \\
We may therefore restrict the resolution of $ \mathbb{C} $ in $ KK^G $ constructed 
in \cite{Voigtbcfo} to obtain a resolution of $ \mathbb{C} $ in $ KK^H $. More precisely, we obtain an exact triangle of the form
$$
\xymatrix{
\res^G_H(C_0(G)) \ar@{->}[r] & \res^G_H(C_0(G)) \ar@{->}[r] & \mathbb{C} \ar@{->}[r] & \Sigma \res^G_H(C_0(G))
}
$$
in $ KK^H $. Recall that the set $ \Irr(G) $ of equivalence classes of irreducible unitary corepresentations 
of $ G $ identifies with $ \frac{1}{2} \mathbb{N}_0 $ and that $ \Irr(H) \subset \Irr(G) $ corresponds to 
the irreducible corepresentations with integral label. Then
$$ 
\res^G_H(C_0(G)) = C_0(H) \oplus C_0^\omega(H) 
$$ 
in $ KK^H $ where $ C_0^\omega(H) \subset C_0(G) $ corresponds to the corepresentations of $ G $ with label in 
$ \frac{1}{2} + \mathbb{N}_0 \subset \frac{1}{2} \mathbb{N}_0 $. It is easy to check that the crossed products $ H \ltimes C_0(H) $ and 
$ H \ltimes C_0^\omega(H) $ are isomorphic to the algebra of compact operators in both cases. 
This holds for both full and reduced crossed products. \\
From these facts it follows that $ H $ is $ K $-amenable, compare \cite{Vergniouxkam}, \cite{Voigtbcfo}. Moreover we obtain an exact sequence 
$$
\xymatrix{
 {\mathbb{Z}^2\;} \ar@{->}[r] \ar@{<-}[d]^\partial &
      K_0(C^*_\max(H)) \ar@{->}[r] &
        0 \ar@{->}[d] \\
   {\mathbb{Z}^2\;} \ar@{<-}[r] &
    {K_1(C^*_\max(H))}  \ar@{<-}[r] &
     {0} \\
}
$$
in which the boundary map can be identified with 
$$
\partial = 
\begin{pmatrix}
n & -n \\
-n & n
\end{pmatrix}
\in M_2(\mathbb{Z}). 
$$
The latter formula is easily verified by inspecting the definition of the resolution considered in \cite{Voigtbcfo}.   
We conclude $ K_1(C^*_\max(H)) \cong \ker(\partial) \cong \mathbb{Z} $ and 
$ K_0(C^*_\max(H)) \cong \coker(\partial) \cong \mathbb{Z} \oplus \mathbb{Z}_n $ as claimed. \qed \\ 
Let us remark that the dual of the quantum group $ \mathbb{F}O(n) $ appearing in the proof of theorem \ref{qautktheory} is 
monoidally equivalent to $ SU_q(2) $ for a certain negative value of $ q $. In our discussion above  
we have restricted attention to $ q \in (0,1] $ for convenience, but the results in \cite{Voigtbcfo} include the case of these 
negative parameters as well.

\section*{Acknowledgments}

This research was partially supported by the Deutsche Forschungsgemeinschaft (SFB 878).

\bibliographystyle{plain}

\bibliography{cvoigt}

\end{document}